\input amstex
\documentstyle{amsppt}
\NoBlackBoxes

\hoffset 1.5cm
\voffset1cm
\topmatter

\def\ttextindent#1{\indent\llap{#1\enspace}\ignorespaces}
\def\itemitemitem
{\par\indent\indent\hskip10pt\hangindent3.5\parindent\ttextindent}

\overfullrule=0pt



\rightheadtext{Traits}
\title TRAITS\endtitle
\author S.~S.~Kutateladze\endauthor
\abstract
Reminiscences about {\it Alexandr Danilovich Alexandrov}
(1912--1999)
\endabstract
\date
\enddate


\address
Sobolev Institute of Mathematics\newline
\indent 4 Koptyug Avenue\newline
\indent Novosibirsk, 630090\newline
\indent RUSSIA
\endaddress

\email
sskut\@ member.ams.org
\endemail

\endtopmatter

\document

Reminiscences and memoirs comprise a special kind of fiction
with  lies and boasts  unavoidable.
The latter were  disgusting for
A.~D.\footnote{Sounds in Russian like ``under'' with
the letters ``n'' and ``r'' omitted and  both syllables stressed equally.}
to an  extent that leaves no room
for envying the authors who provide their written recollections about A.~D.

It happened so that I had a privilege and honor of
constant communication with A.~D. from  the end of the 1970s
up to his death. Writing reminiscences is by far much easier
after  many years' elapsing. However, my elder friends
had managed  to convince me to reflect
some details of the Siberian period of A.~D.'s life.

I have had many opportunities of writing about A.~D.
in the traditional (and not fully traditional) forms of scientific
publicism. I am happy that he   never  reproached me for this,
and so I guess that  I may skip surveying  his
scientific  contribution.

Many events in which I observed A.~D. and sometimes participated in
secondary roles
were not so long ago as to become an impartial history.
Not all of them deserve inspecting over for revival and plunging
into once again.

Perusing my personal archives, I decided to select just a few items
reflecting those traits of A.~D.'s personality
that were revealed in our contacts.

I will be glad if the lessons of A.~D.'s life help someone to hold
on or to settle some pending  crisis as they have readily done for me...

\head
Anger and Self-Criticism
\endhead

A particular trait of A.~D.
I wish to emphasize is the
physiological reaction of anger
to danger, assault, or offence.
It is well known that
These circumstances  are well known to
bring about the emotion of fear (pale face, cold wet, etc.)
The military commanders of the ancient times often enrolled in their
forces
the warriors 
whose  reaction to danger  was anger.

A.~D. exhibited the classical examples of the emotion of anger:
his face reddened, the chest threw out, and he showed the bared teeth.
A.~D.  understood quite perfectly how he intimidated those who
provoked his anger.
At that I never saw any unjustified fits of his anger.
Many years of acquaintance with A.~D. cultivated the strong stereotype:
Everyone hating A.~D. is a potential if not complete scoundrel.
As regards his students, friends, and relatives, A.~D.
was exceptionally kind, even tender, very attentive, and scrupulous.

A man of passion, A.~D. always remained self-critical.
I had an opportunity to write that self-criticism is a necessary
test for intelligence.
Every now and then A.~D.  reconsidered his attitude to
people and events in accord with the ideals of morality
he proclaimed: universal humanism, responsibility, and
scientific outlook.

As a small illustration, I can recall that A.~D.
voted against the admittance of my Ph. D. thesis
to  the formal procedure of public maintenance
in 1969. Moreover, he supplied no motivation whatsoever.
An open negative vote of an academician
happens rarely on such a trifle occasion as the admittance
of somebody's Ph.~D. thesis.
My thesis was submitted in analysis
under the  title {\it Related Problems of Geometry and Mathematical Programming}.
Its topic was close to the research of
A.~D. Alexandrov in the theory of mixed volumes
and the research of L.~V. Kantorovich
in optimization and ordered vector spaces.
Clearly, I was not the only person impressed by the unmotivated demarche
of A.~D.

I thought that my article could be of interest
to A.~D.
(the formal review of a ``leading mathematical organization''
was written by V.~A.~Zalgaller; and my main technical
result was an extension of one unpublished idea by Yu.~G. Reshetnyak
in measure theory).  I was rather nervous making my talk at the public
maintenance. Using an overhead projector in a semi-darkened hall,  I cast
a casual glance towards A.~D.
When I had told that my thesis
uses the theory of mathematical programming by
L.~V.~ Kantorovich and the theory of surface area measures by
A.~D.~ Alexandrov, there was some noise from the side benches:
A.~D. rose and strode out. It is easy to imagine
how confused I was after that. However, the vote was unanimous.

After many years, when we had been close with  A.~D.
for a long time, I reminded him of this story.
He rebuffed immediately:
``This never happened at all.'' (You should know  A.~D.
to understand his answer
properly: when he forgot or doubted something, he
always said: ``Don't remember.''   Replying in other words,
A.~D. had declared the whole episode nil and void.)

It is a pleasure to recall that I had received satisfaction from
A.~D. in due course. As a result of some bizarre machinations of
the All-Union Attestation Committee in the 1970s, my
Sc.D. thesis was sent to extra referral
despite its formal approval at the corresponding section of the
Committee on the recommendation by
E.~M. Nikishin. Happily, the appointed ``black'' opponent
was A.~D., and I received his appraisal for the isoperimetric
problems with arbitrary constraints on mixed volumes.

\head
M.~A. Lavrent\kern1pt$'$\kern-1.5ptev and a Book on the Methodology of Mathematics
\endhead

Narrating about his participation in the
ideological battles of the 1940s and 1950s,
A.~D. always spoke about the tactics of preemptive
blows. One of them deserves recalling.

The Academy of Sciences of the USSR had printed in 1953
a huge volume of about thousand pages under the title:
{\it Mathematics: Its Content, Methods and Meaning.}
The Editorial Board of the volume comprised
A.~D.~Alexandrov,
A.~N.~Kolmogorov, and M.~A.~Lavrent$'$ev.

The eighteen chapters of the book were intended to the general public
and written by thirteen authors. The list of the latter
contained  I.~M.~Gelfand, M.~V.~Keldysh,
M.~A. Lavrent$'$ev, A.~I.~Mal$'$tsev,
S.~M.~Nikol$'$ski\u\i{}, I.~G.~Petrovski\u\i{}, and S.~L.~Sobolev.
The  run of 350 copies was exceptionally small
those days. Besides, each copy was enumerated and the front page contained
the index of the copy
and the extraordinary signature stamp ``Published for Discussion.''

Sufficiently many copies of this book
were printed free of classification
only in 1956, and the book became an issue
in the world mathematical literature.
Suffice it to say that the translation of this book was  reprinted
thrice in the USA (the last time in 1999).

Clearly, such an extraordinary collection
had rather nontrivial reasons for its compilation.
The aim of this project consisted in defending
mathematics from the antiscientific attacks that were typical of those days
in the Soviet Union.

To strike a severe preemptive blow on the pseudo-scientists
of marxism which try to harass the development of
science in this country and to get rid of them
possibly for ever  was an almost successful plot of
the book.
Avoiding strictly professional nuances,
the world-renowned leaders of mathematics gave  in this book
a detailed and thorough analysis of such fundamental
general aspects of their science
as the subject of mathematics and the nature of mathematical abstractions,
interaction     between pure and applied mathematics,
relationship between mathematical research and practice, etc.
The book remains to stand as one of the heights
of the methodology of mathematics.

The soul of this project was A.~D.
In addition to the two special chapters on curves and surfaces and
on abstract spaces, he  made a ``promising beginning'' ---
wrote the lengthy introductory   chapter
``A General View of Mathematics'' with an impressive analysis
of the challenging philosophical problems of mathematics.

The work on this book drove A.~D. and M.~A.~Lavrent$'$ev closer.
By the invitation of M.~A.~Lavrent$'$ev,
A.~D.~Alexandrov joined the staff of the
Siberian Division of the Academy of Sciences of the USSR
in 1964.
A.~D. was proud of the fact that M.~A.~Lavrent$'$ev
had solely nominated  him as a candidate to a full member
of the Academy
and freed him from all bureaucratic formalities.
When A.~D. became aware that L.~V.~Kantorovich
was nominated for the same vacancy, he began
to refuse to participate in the elections. However,
M.~A. managed to convince A.~D. to stop refusing.
Sage Mikhail Alekseevich turned out to be right: Both were happily
elected (the Bylaws of the Academy made room for such an outcome
those days).

\head
Bertrand Russell and a Preventive War Against Russia
\endhead

At the end of the 1970s the plan was under discussion of publishing a volume of the articles of A.~D.
on the general problems of science and other articles of publicism.
This plan  led finally to his book {\it Problems of Science and a~Scientist's Standpoint}.
The release candidate No.~1 was surely the article ``A General View of Mathematics.''
A.~D. asked me to look it through for shortening.
Reading the article thoroughly, I felt much doubt about the
following excerpt:
\medskip
\itemitem{}{\eightpoint\sl\indent
In the bourgeois society we can encounter the scientists
that turn into obscurantists professing political reaction and
antiscientific obfuscation rather than progress  and knowledge.
An example of this degeneration is one of the founders of the so-called
``logical positivism''---Russell, an English philosopher and mathematician.
He declared fifty years ago that
``Mathematics may be defined as the subject where we never know what
we are talking about, nor whether what we are saying is true.''
In other words, mathematics has no real content according to Russell.
The real content of his own views Russell revealed completely
when  he began to call for  atomic war against the Soviet Union
several years ago.  A perverter of science and self-conceited
epigone of forlorn idealistic systems who instigates
mass destruction---that is the true face of this ``logical positivist.''}

\medskip\noindent
In my opinion of those days Russell was
one of the leaders of the Pugwash movement, a dedicated peace warrior, and a Nobel prize
winner. In no way  he was reminiscent of a perverter of science who
instigates mass destruction.
Frankly speaking, I thought that A. D. swallowed a~tasty bait
of the propagandists of the CPSU in the first years
of the Cold War.

With a hardly concealed spite I told  A.~D. that the reader
needs a precise reference to the words of Russell
and smugly requested his explanations.
In fact, I attacked him impudently in the trite style
of the ``presumption of dishonesty.'' He was definitely offended.
He snapped back sharply that the episode did take place but
he could not remember any details. I must confess that
these explanations convinced me of nothing at all.

In the new millennium I tried to use the omnipotence of the Internet
to settle the problem finally by  search machines.
Without any effort, I found out that one of Russell's
statements about the A-bomb appears in textbooks as a standard
example of  a~``false dilemma.''

\medskip
\itemitem{}{\eightpoint\sl\indent
Either we must have war against Russia before she has
the atom bomb, or we will have to lie down and let them govern us.}
\medskip\noindent
Also, Russell included in his book
{\it The Future of Science, and Self-Portrait of the Author\/}
published in~1959 the following interview he gave for BBC Radio:

\medskip
\itemitem{}{\eightpoint\sl
Q. Is it true or untrue that in recent years you advocated that a
preventive war might be made against communism, against Soviet
Russia?"}
\itemitem{}{\eightpoint\sl
RUSSELL: It's entirely true, and I don't repent of it now. It was not
inconsistent with what I think now.... There was a time, just after
the last war, when the Americans had a monopoly of nuclear weapons and
offered to internationalize nuclear weapons by the Baruch proposal,
and I thought this an extremely generous proposal on their part, one
which it would be very desirable that the world should accept; not
that I advocated a nuclear war, but I did think that great pressure
should be put upon Russia to accept the Baruch proposal, and I did
think that if they continued to refuse it it might be necessary
actually to go to war. At that time nuclear weapons existed only on
one side, and therefore the odds were the Russians would have given
way. I thought they would ... . }
\itemitem{}{\eightpoint\sl
Q. Suppose they hadn't given way.}
\itemitem{}{\eightpoint\sl
RUSSELL: I thought and hoped that the Russians would give way, but of
course you can't threaten unless you're prepared to have your bluff
called.}
\medskip\noindent
It is a pity that A.~D. will never hear the words of my repentance.

\head
A.~P. Aleksandrov and a Polemic about an Article by
N.~P.~Dubinin
\endhead

A.~D. was engaged in defense of science and particular scientists
in his Siberian period as well.
Many persons he drew out of the screw presses
of the scientific and would-be-scientific  rascals
who made their careers in the 1970s and 1980s.
I am reluctant to tell these stories
whose analogs are familiar to the majority of
scientific groups in this country.

What I want to recall here is the valiant standpoint of
A.~D. in regard to the article by N.~P.~Dubinin
``Biological and Social Heredity'' which was
published in
{\it The Communist}~ (1980:{\bf11}).\footnote{The
official journal of the Central Committee
of the Communist Party of the USSR.~({\bf S.~K.})}

A.~D. had appraised this composition as
an ``outstanding piece of antiscientific literature.''
I am convinced that to read the article by N.~P.~Dubinin and
the relevant controversy
is as vital for a young scientist of any specialty as the
perusing
of the shorthand record of the notorious
August Session of the Lenin All-Union Academy of  Agricultural Sciences 
(VASKhNIL in the Russian abbreviation) which took place in
1948.

Avoiding to  narrate the whole composition of
N.~P.~Dubinin, I just pinpoint  one of the
ideological conclusions of his article:

\medskip
\itemitem{}{\eightpoint\sl\indent
Without clear understanding of the genuine
scientific basis for the
problem of  man, it is impossible
to properly place the vicious essence of neoeugenical
ideas in a disguise of  new discoveries in
natural sciences and in particular in molecular biology
and genetics.  Moreover, this problem is such that
the coincidence of the truth criterion and the party
spirit is most conspicuous here.}
\medskip\noindent
A.~D. found  primarily repulsive the attempt at
making the party spirit the test for truth
and refused to consider it as a slip of the tongue.
His worst premonitions came true:
the editorial comment on the discussion around the article
by N.~P.~Dubinin  read later in  {\it The Communist} (1983:{\bf14}):

\medskip
\itemitem{}{\eightpoint\sl\indent
The main criterion for evaluating the
philosophical meaning of pieces
of theoretical research is their ideological orientation,
the purity of the  class-charact\-eristic, ideological, and methodological
standpoints.}
\medskip
\noindent
Practice as the ultimate test for truth was doomed for a pompous
funeral and complete oblivion.

A.~D. tried to  profess his views of
the article by N.~P.~Dubinin actively: he made talks on methodological
seminars in various scientific institutions
and attempted in vain to publish his arguments.
Fortunately (this happened quite rarely to A.~D.),
he was supported by
A.~P.~Aleksandrov who then held
the position of the President of the Academy of Sciences
of the USSR
and let A.~D. get the floor at the General Assembly of the
Academy of Sciences of the USSR on November~21, 1980
(a version of the speech  of A.~D. and the reply by
N.~P.~Dubinin are published in the  {\it Herald of the
Academy of Sciences of the USSR} (1981:{\bf6})).

A.~D. always emphasized that the cause of science is to find out
``how the thingummy's actually going on.'' He pursued the same approach
in this particular case:
\medskip
\itemitem{}{\eightpoint\sl\indent
The real problem consists in
studying which sides of psyche depend on heredity or
social environment and to which extend.
However, N.~P.~Dubinin closes this problem as regards
normal persons, leaving it open only for
medical genetics in regard to abnormal persons.}
\medskip\noindent
A.~D. told me after the Assembly that Anatoli\u\i{} Petrovich
answered to A.~D.'s application for having the floor
as follows:
``Do you want to bite off Dubinin's head right away or after
the break?''
As far as I could remember, A.~D. was eager to accomplish
the task immediately...
These days  A.~D. gave me a galley proof of the draft of his speech.
Below I present  the end of this manuscript which remained unpublished by now:
\medskip
\itemitem{}{\eightpoint\sl
I have said now what I wished to say,
but I harbor heavy doubts:
maybe, it was unnecessary to speak all this out and
in so strong words at that.
In fact, it is clear that the attempts of Academician
Dubinin will not affect serious scientists and practitioners.
Therefore, they will hardly influence our biology
and medicine.}
\itemitem{}{\eightpoint\sl
\indent However, this view is not fully
accurate.
Academician Dubinin used a~high rostrum
and it is not completely excluded that
some assistant professor reading human genetics  in
some medical college will be called to
responsibility
for ``attempting''---in the words of
Dubinin---``to revise and waste the marxist teaching
of the unique social nature of  man.''}
\itemitem{}{\eightpoint\sl
\indent Furthermore, the question is  placed on the agenda of
the honor of science and  our personal honor:
Do we agree to yield to the resurrection of the vicious
style and battle against science which reigned here
about thirty years ago?}
\itemitem{}{\eightpoint\sl
\indent Everyone can make a mistake and even
speak up in rigmarole.
 What really matters in the long run are
 the basic principles of science themselves
rather than some particular mistakes.
Marx observed that anyone  I call a base person who strives to
adjust science to  external and alien aims%
---irrespective of whatever delusions   science might cling to.}
\itemitem{}{\eightpoint\sl
\indent The objects of concern are precisely  the main principles
of scientific research:   impartiality and honesty.
We cannot let them  be dismissed so
loudly and impudently.}

\head
S.~L.~Sobolev and a Polemic about  an Article by L.~S.~Pontryagin
\endhead

The year 1980 was rich in events!

The journal {\it The Communist} (September 1980:{\bf14})
published the article by  L.~S. Pontryagin
``About Mathematics and the Quality of Teaching Mathematics.''
This composition still arouses  the emotions as sharp
as those stirred up by the article of N.~P.~Dubinin.
Moreover, both in the same volume of the journal produce
an unforgettable adore.

The article by L.~S.~Pontryagin was supplied with
a routine editorial comment that explained the genuine meaning of
the article to those who
failed or hoped to fail to grasp it:

\medskip
\itemitem{}{\eightpoint\sl
...the author is right in opposing vehemently
not only the exceeding devotion to abstract constructions
in teaching mathematics and in
mathematics itself but also
the pseudoscientific speculations related to the false
treatment of its subject.}
\itemitem{}{\eightpoint\sl
\indent Noncritical adoption of foreign
achievements in relatively new
branches of mathematics and hypertrophy of  general
importance of these achievements to science as a whole
have led to overrating the results of many
mathematical studies
and in some cases to  the idealistic
treatment of the essence of the subject of this science,
to the absolutization of abstract constructions, and
to the belittling of the gnosiological role of practice.
Exceeding devotion to  abstractions
of the set-theoretic stance
has started disorienting the creative
interests of students and academic youth.}
\medskip\noindent
It was impossible to consider such a rhetoric
casual and innocent. Indeed,
{\it The Communist}~{\bf18} had published  a note by
Academician I.~M.~Vinogradov, Director of the Steklov
Mathematical Institute of the Academy of Sciences of the USSR.
This note said in particular that

\medskip
\itemitem{}{\eightpoint\sl
The Scientific Council of the Steklov Institute
is satisfied with the statement of the journal
{\it The Communist\/} in the form of the article by Academician
L.~S.~Pontryagin... The Scientific  Council of the Steklov Institute
supports the statement of the journal
{\it The Communist\/} and believes that it will serve
the cause of
improvement of teaching in secondary school...}
\medskip\noindent
I find it appropriate to cite a few lines from my diary
for reconstructing the intensive but stale atmosphere of that
span of time.

\medskip
\itemitem{14.10.}  A.~D.  called me in the evening and told
about an article in {\it The Communist\/}
{\bf14}: Pontryagin vs. secondary school,
S.~L. + L.~V. [Kantorovich] + an editorial comment on
idealism in mathematics.

\itemitem{15.10.}  M.~A.~Lavrent$'$ev passed away.

\itemitem{18.10.}I read the nasty article by Pontryagin in the morning
and dropped in on A.~D. in the evening to talk this over.

\itemitem{24.10.} The ninth day---the funeral of Mikhail Alekseevich...

\itemitem{25.10.} Zelmanov was flunked by secret ballot in the Institute.
We discussed this and  {\it The Communist\/}  with  A.~D. at length.

\itemitem{26.10.}  A.~D. dropped in on me and gave me the second part
of his textbook. Then I dropped in on him. A.~D. wants to retire.

\itemitem{30.10.} ...G.~P. [Akilov] crashed his car but slightly.
A.~D. told me that S.~L. has written a reply to  {\it The Communist}.
Yu.~F. [Borisov] called me about extreme points.

\itemitem{3.11.} At V.~L. [Makarov]'s seminar in the morning.
Then I visited  S.~L.---about  {\it The Communist}.
He showed me  his reply.
Next---at  S.~L.'s seminar with a Vietnamese.
Then another conversation with S.~L. about the article (in the breaks
I talked to A.~D.). S.~L. spoke eloquently but slightly incoherently about sets and cardinality.

\itemitem{4.11.} S.~L. called me about arranging a meeting of the
Scientific  Council
vs. Pontryagin.

\itemitem{5.11.} [The Scientific  Council] unanimously condemned the
Thesis Maintenance Co\-uncil. The speakers were
Serebryakov, Yu.~G. [Reshetnyak], A.~D., et al....

\itemitem{10.11.} The whole day was full of discussions with S.~L. and A.~D.
about Pontryagin and also
about Reshetnyak and Zelenyak (in view of a~ scandalous meeting of the
Academic Council
to take place tomorrow [in NSU]).

\itemitem{11.11.} An anniversary of Bourbaki. Yu.~G. was cancelled by
 40\% ...

\itemitem{12.11.}  A seminar with  A.~D. about Lenin's speech at
the  III Convent + vs. Pontryagin... Yu.~G. discussed $L_p$  with
me.

\itemitem{24.11.} I looked for the list of the members of the Council with
L.~M. [Krapchan]. A.~D. has arrived---he spoke
against  Dubinin at the General Assembly. Dubinin replied...
A seminar about attractors with Ustinov (from Obninsk).

\itemitem{28.11.} [The Scientific  Council] unanimously supported the appeal by Zelmanov.
Celebration of the 20 years of the M[athemati\-cal-]E[conomical]D[epartment].

\itemitem{3.12.} I dropped in on S.~L. about the resolution.
He told that he will move it himself.

\itemitem{8.12.} I dropped in on
S.~L. with A.~D., V.~A. [Toponogov], and
V.~V. [Ivanov].  Discussion  vs. Pontryagin.

\itemitem{12.12.} [We were] pretty close to adoption of the anti-Pontryagin resolution
[at the philosophical-methodological seminar]. I dropped in on A.~D.
in the evening to  talk this over.

\itemitem{15.12.} I discussed the resolution  with S.~L.
Then at his seminar...
Bokut$'$ called me in a break about his troubles. Shirshov
recommended Ershov for the party membership...

\itemitem{23.12.} Talking everything over the whole day out with Tikhomirov
who just arrived. Mainly in the anti-Pontryagin mood.

\itemitem{24.12.} Thesis maintenance: [V.~N.]~Dyatlov 18-0=0 and
[G.~G.]~Magaril
[-Il$'$yaev] 17-0=0...
Everything was pretty nice... Booze\&noise at Dyatlovs'...

\itemitem{25.12.} Short discussion with A.~D. in the morning...
[The Scientific  Council] adopted the anti-Pontryagin text + there will be a letter to
{\it The Communist\/} to be prepared by
A.~D. +  Yu.~L. [Ershov] + [S.~I.] Fadeev!
\medskip

\noindent
Such were the circumstances we lived in those days.

I remember the extraordinary stamina of A.~D. (which was predictable)
and Serge\u\i{} L$'$vovich (which was unexpected for me).
The latter startled me on November~3, giving his reply to
{\it The Communist}:
``I am interested in your opinion but you should bear in mind that I have
already mailed my reply.''
On the same occasion he showed me a copy of an analogous letter
to somebody in the leadership of the Central Committee of the CPSU
(seemingly, this was M.~V.~Zimyanin).

Many participants of these events are still alive.
Some of them have changed for the better (and the rest of them
still have a good chance to do the same).
That is why I am reluctant to describe all details of the
vehement struggle for a noble answer to the article by
L.~S.~Pontryagin.
I mention only that the crucial ingredient was
the titanic joint efforts of Aleksandr Danilovich
and Serge\u\i{} L$'$vovich.

In result, the Scientific  Council of the Institute of Mathematics
unanimously (sic!) adopted at its meeting of December~25, 1980 the
resolution that read in particular as follows:

\medskip
\itemitem{}{\eightpoint\sl\indent
The Scientific  Council announces its disagreement with those
who informed the Editorial Board of  {\it The Communist\/}
about the situation in the science of mathematics
which gave grounds for  the editorial comment on the
article by  L.~S.~Pontryagin  to make accusations of
the noncritical adoption
of foreign
achievements, formalistic craze, disorientation
of academic youth,  and  the false treatment
of the subject of mathematics.
Mathematics is a unique whole and
deterioration of its fundamental more abstract part
resembles  proscription of
chromosome heredity theory,
treatment of cybernetics as a ``science of obscurantists,'' and
prohibition of using mathematical methods in economics on the basis of false
pseudoscientific
arguments.
Mathematics is a very serious matter of paramount importance
for the development of our society.
Therefore, treating it and judging it requires
great responsibility.}
\medskip\noindent
There was some cool in the relations of
A.~D. and S.~L. that year (but I am disinclined to reveal the reasons behind this yet).
Therefore, it happened so that
the drafts of the resolutions were prepared with me acting as
an intermediary.
I keep these drafts with the scars of those ``shuttle operations'' in remembrance of the unforgettable material lesson of
struggling for scientific truth.

It is worth  observing that  E.~I. Zelmanov
whose Ph.~D. thesis was rejected by secret vote as
mentioned above  acquired the Fields Medal a few years later.

The standpoint of Serge\u\i{} L$'$vovich was reflected by
{\it The Communist\/} in the phrase: ``Comments are still coming.
Among them some are  written in a polemic  style:
the letters by Academician
S.~L.~Sobolev,
Assistant Professor P.~V.~Stratilatov, and
Professor Yu.~A.~Petrov.''
The chant
``Academician Sobolev,
Assistant Professor Stratilatov, and Professor Petrov''
was our catch-phrase
for a few years.

We attempted to print  a booklet with the resolution of the
Scientific  Council and a detailed version of the report
by A.~D.~Alexandrov ``About the Article by L.~S.~Pontryagin
in  {\it The Communist} (1980:{\bf14}).''
Our attempts were unsuccess\-ful---we were opposed by
V.~A.~Koptyug.
\footnote{The Chairman of the Presidium
of the Siberian Division (1980--1997). ({\bf S.~K.})}

A.~D. showed me a personal memo by
V.~A.~Koptyug in which the latter---a censor (sic!)---reproached  A.~D.
for a~``persecutor's tone'' and refused to publish
the report.

Despite this the scientific community became aware of
the standpoint of Siberian mathematicians:
at Sobolev's request the copies of the resolution and
A.~D.'s report were  sent to the principal
mathematical institutions of this country.

Something similar happened later to  A.~D.'s book
{\it Problems of Science and a~Scientist's Standpoint\/} whose publication
was procrastinated by the chiefs of the Siberian Division
and became possible only after interference of P.~N.~Fedoseev
who knew  A.~D.
rather well  and strictly obeyed
academic etiquette in this matter.

\head
Sic Transit...\\
or  Heroes, Villains, and Rights of Memory
\endhead

April 25,  2003  is the date of the centenary of the birth of
Andre\u\i{} Nikolaevich Kolmogorov.
The personality and creative contribution
of this genius man to the world science and Russian
culture are so eminent that the tiniest bits of
recollections of anything related to him
might be of avail to those  pondering over life and its principles.

For many years I have heard requests of my friends and colleagues
to present for the public my whatever partial overview of
the circumstances and events invoked
by Merzlyakov's  article ``The Right of Memory''
and in  particular the  polemic between A.~D. Alexandrov
and L.~S. Pontryagin  this article had stirred up.
The story to tell is rather ugly and to plunge into it again,
reviving  the bygones, brings about much discontent and displeasure.

Unfortunately, the historical nihilism of these days
intertwines rather  tightly with  nihilism in morality.
``The past crimes are buried in the past. The past is absent at
present. Therefore, the past crimes are absent now. So,
let bygones be bygones.'' This sophism brings about the opinion  that
nobody could recall  and take into account the crimes of the past
in view of  the period of limitations.
This is correct but partly. The murderer remains a murderer for ever
irrespective of whether or not he committed a negligent homicide
and was relieved from persecution or  has served his punishment
and lives with no record of conviction.
The thief is still a thief although she returned back
the things she had pilfered
and was relieved from punishment.
No fact of assassination or theft
is ever repealed by whatever decisions about it.
No error disappears unless it has been repaired.
Always evil is to forget the past and its lessons...
These arguments drove me to the decision of
narrating about this gloomy episode of the past.

Merzlyakov's article  appeared on February 17, 1983 in the newspaper
{\it Science in Siberia\/} of the Presidium of the Siberian Division
of the Academy of Sciences of the USSR.
Yu.~I. Merzlyakov (1940--1995), an established algebraist, a Sc.D. and professor,
had a bit of reputation in the theory of rational groups.
He was not an ordinary personality
devoid of literary and other gifts and so won quite a few admirers.
His article served many years as a credo of the Novosibirsk
branch of the notorious ``Memory'' society, an informal nationalistic
group sprang to life
in the early years of Gorbi's perestro\u\i{}ka.

To grasp the undercurrents of Merzlyakov's article completely
is practically impossible for anyone far from the
Russian mathematical life of those days.
Moreover, the understanding of and attitude to this text varied
drastically from capitals to province.
Despite this, all Russian mathematicians  clearly saw the implication of the
following excerpt of Merzlyakov's article:

\medskip
\itemitem{}{\eightpoint\sl
Academician Lev Sem\"enovich Pontryagin is a brilliant
exemplar of a scientist and patriot of these days.
The International Astronautical Federation
elected him an honorary member side by side  with the
cosmonauts Gagarin and Tereshkova for his outstanding scientific contributions.
Skipping any description of all aspects of the versatile
activities of L.~S. Pontryagin,\footnote{The initials seem abundant
to the English eye
but they reflect the style of the Russian polemics in which the presence of
initials brings about some extra  respect to the persons in question whereas
the absence of initials clearly demonstrates slight indifference,
disrespect, or even neglect. Every Russian professor knows that the initials
of Gagarin are Yu.~A., and the initials of Tereshkova are V.~V.
To keep the flavor of the polemic I preserve the authors'
rules for placing initials in the Russian originals throughout.~({\bf S.K.}) }
I will dwell upon a single
problem of a national-wide scale, the teaching of mathematics in
secondary school. It was exactly L.~S. Pontryagin
who vehemently pointed out, in particular on the pages
of {\it The Communist},
\footnote{{\it The Communist}, 1980:{\bf14}, p.~99--112}
the  evil implications of the sharp turn to the course of the extreme formalization
of mathematics which  was imposed on our schools in 1967 and
oriented consciously or unconsciously
to the accelerated intellectual
development (with an equally fast achievement of
the utmost limits of creativity) nontypical
of the majority of the country's population.
The flood of responses to
the statement by L.~S. Pontryagin
\footnote
{{\it The Communist\/}, 1980:{\bf18}, p.~119--121;
1982:{\bf2}, p.~125--126} has demonstrated that
his criticism was quite timely and fair.
In particular, Vice-President of the Academy of Sciences of the USSR
Academician A.~A. Logunov ascertained on the
session of the Supreme Soviet of the USSR in October 1980
that there is a grieve situation about the
teaching of mathematics
in secondary school and to learn mathematics from the present-day
textbooks ``can destroy any interest  in not only
mathematics but also  exact sciences in general.''
(I remark parenthetically that the leader of the reform
received the prize of 100,000 dollar in 1980 from
the state with which the USSR had severed all diplomatic relations
exactly in the year of the beginning of the reform.)}
\footnote{{\it Notices of the AMS\/}\,(1981)\, {\bf 28}:1, p.~84}
\medskip
\noindent
The rest of the article was mainly inspired
by the outright scandalous situation in the midst of
logicians and algebraists of Novosibirsk
and in the whole mathematical community of Siberia either.
The point was that the retirement of S.~L. Sobolev was pending  from the
position of the director of the Institute of Mathematics of the Siberian
Division of the Academy of Sciences of the USSR. This evoked
the battles for power and better places under the sun which
were typical of the academic community of those days.

I am disinclined to dwell upon the other details of
Merzlyakov's article since I fully agree with
the estimate of Sobolev who expressed his attitude to
the hysterics by Merzlyakov as follows: ``The role of Savonarola
befits no twentieth-century scientist.''

Sobolev forwarded his sagacious and valiant letter from Moscow
to the management of the Institute
on March~9.  He rejected the slander against Kolmogorov
and justly gave a negative estimation of the whole article.
I had an opportunity to read this hand-written page of a copy-book
which unfortunately was unwelcome by some of the addressees,
concealed for a long time, and made public by S.~K. Godunov only
after fierce battles and conflicts
at the meeting of the Academic
Council of the Institute  on April~18.
The principled and uncompromising position of Sobolev
seemed to the many less important than the opinion of local party
leadership. A few iterations under the pressure of petty  communist bonzes
brought about the official position of the management of the
Institute which recalled the merits of Kolmogorov while
observing that Merzlyakov appropriately posed the problems of patriotism.

Patriotism and slander... A notorious mixture...

Some unpleasant general thoughts are in order now
about professionalism and mathematicians.
Professionalism requires absolute devotion to
profession and, absorbing personality, tends to
impoverish the latter.
Professionalism appears  amidst mathematicians rather early
whereas the upbringing of necessary moral qualities
is often far from a fast and easy matter
(mathematicians are next of kin to sportsmen in this respect).
Of little secret are the elements of gossip,
jealousy, and envy  encountered the world over even among the first mathematicians.
Hatred to the gifts of the others is often mixed
or replaced with  xenophobia, racism, antisemitism, and  similar elements
of the same sort. These phenomena are still  far from rare
nowadays. The oversensitive reaction to the slightest traits
of the presence or absence of antisemitism
was and still is a litmus test of ``friend-enemy'' in
Russia irrespective of whether this  is right or wrong.
I believe that to grasp correctly the tension of the events
after Merzlyakov's article is impossible
without the clear understanding of the above circumstances
of the Russian life.

By the way, somebody told me that
the then editor-in-chief of the newspaper  {\it Science in Siberia\/}
tried to justify  himself on explaining that he had slightly
deviated from the standard routine of
accepting materials for publication in order to insert
Merzlyakov's article in the issue on the Day of the Soviet Army
because he viewed it as exceptionally patriotic.
In our midst we have called these views ``slanderous patriotism'' since then.
Mixing  love for the Fatherland with slander is
always characteristic of ``the last resort of a scoundrel.''

The Moscow mathematical community reacted  to
Merzlyakov's article immediately and adequately in general.
The understanding prevailed that the lampoon
could strike the health of Kolmogorov which
was already shaken seriously. Surely, nobody showed the newspaper to
Andre\u\i{} Nikolaevich  but his 80th anniversary approached rapidly
and Merzlyakov's article could provoke
some undesirable predicaments: for instance, there might have been no
ceremonial decoration from the government
which could be noticed by Kolmogorov, stirred up his
analytical interest and investigation with
possibly unfavorable aftereffects to his health.

Another circumstance helped to the spreading of a noble reaction:
The article appeared on the eve of the General Assembly
of the Academy of Sciences of the USSR in Moscow
where several copies of the issue of the newspaper were
delivered immediately.
The exceptionally sharp reaction against
slander and the style of a political informant
was revealed by the leading mathematicians: A.~D. Alexandrov,
S.~M. Nikol$'$ski\u\i, S.~P. Novikov, Yu.~V. Prokhorov, S.~L. Sobolev,
L.~D. Faddeev, and many others.

Already on March 14 there appeared the first written response by
Alexandrov with an analysis of Merzlyakov's article. Characterizing the
article as objectively anti-Soviet and subjectively base, Alexandrov
demonstrated  the necessity of terminating all instances of slander and
political insinuation. Closing his response, Alexandrov wrote:
\medskip
\itemitem{}{\eightpoint\sl\obeyspaces\indent
Yu.~Merzlyakov himself has clearly deserved the
right of memory. Since at least some of his
statements are so evil and monstrous that  might
go down in history....}
\itemitem{}{\eightpoint\sl\obeyspaces\indent
We have thus seen that Merzlyakov's article
is an objectively anti-Soviet, subjectively base,
rude, and antipatriotic composition,  its every appeal to patriotism notwithstanding.}
\itemitem{}{\eightpoint\sl\obeyspaces\indent
Let us abstain from judging the author severely
but rather pity him since we observe
an indubitable pathological case.
Only a perverted mind and turbid imagination
can bring about such a flood of insolence and mud!
Renegades, domestic emigrants, immature  moral viewpoints halfway from
amoeba to cave-dweller, a shitting bull, a beast, a toady-like mediocrity
of a petty shop-keeper and, to crown all these,
the monstrous image of villains that crawl to loot the wounded
as description of the ``horde'' of  scientific workers
and, in particular, his fellow colleagues.
Well, that is the limit: an obvious  pathology.}
\medskip
\noindent
We are to pay tribute to the Mathematics Division
of the Academy of Sciences of the USSR  and personally to
Yu.~V. Prokhorov who was an initiator and editor of the following
Resolution of the Bureau of the Mathematics Division as of
March 25, 1983:.

\medskip
\itemitem{}{\obeyspaces\eightpoint\sl
\indent Academician Yu.~V. Prokhorov informed the body
about a recent article in the weekly newspaper
{\it Science in Siberia\/} of the Presidium of the SDAS\footnote{
The abbreviation of `` Siberian Division of the Academy
of Sciences.'' ({\bf S.~K.})} of the USSR (No. 7 of February 17, 1983)
by  Yu.~I. Merzlyakov, Sc.D.  on the staff of the Institute
of Mathematics of the SDAS of the USSR. This article contains
a uniquely decipherable insinuation against Academician A.~N. Kolmogorov,
 an outstanding Soviet scientist.}

\itemitem{}{\obeyspaces\eightpoint\sl
\indent The floor for discussion was taken by
Academicians S.~M. Nikol$'$ski\u\i, V.~S. Vladimirov,
S.~P. Novikov, A.~A. Samarski\u\i, S.~L. Sobolev, and L.~D. Faddeev;
Corresponding Members of the Academy of Sciences of the USSR
A.~V. Bitsadze, I.~M. Gelfand, A.~A. Gonchar, and
S.~V. Yablonski\u\i. All of them unanimously condemned
the indecent insinuations  of Merzlyakov's article
and qualified them as slander against one~outstanding scientist
and patriot. It was also observed that the article contains
insinuations against a number of other Soviet mathematicians.}

\itemitem{}{\obeyspaces\eightpoint\sl
\indent The Bureau of the Mathematics Division of
the Academy of Sciences of the USSR
HAS DECIDED}

\itemitemitem{1:}{\eightpoint\sl
to observe that the article of Sc.D. Yu.~I. Merzlyakov
``The Right of Memory''  in the newspaper
{\it Science in Siberia\/} of the Presidium of the SDAS of the USSR
contains slander against one outstanding scientist/mathematician and
Soviet patriot;}

\itemitemitem{}{\eightpoint\sl
to observe that the article contains  a~number of
indecent attacking allusions to
other Soviet mathematicians.}

\itemitemitem{2:}{\eightpoint\sl
to call upon the Presidium of the Siberian Division of the Academy
of Sciences of the USSR to take due measures pertinent to Item~1.}

\itemitem{}{\eightpoint\sl
The Resolution of the Bureau of the Mathematics Division of
the Academy of Sciences of the USSR was adopted by a~unanimous
vote.}

\medskip
\noindent
The bushes of provincialism
were already full-fruited in Siberia those days, and the solicitude
for the honor, dignity, and health of Kolmogorov together with
counteraction against the filthy things like antisemitism
seemed to the chosen few to be negligible as compared
with the prevailing sentiments for their own career, success, fame, and prosperity.
The following story of Alexandrov looks like a joke nowadays:
one of the top bosses of the Siberian Division responded to the protest
and indignation against Merzlyakov's article with
the sincere question: ``Who is that Kolmogorov guy?''
One can easily imagine our reaction...

On March 28 there was a meeting of the Presidium of the
SDAS of the USSR.  The official letter of the Institute,
bearing the signatures of the three deputy directors and the party
secretary, was announced together with the second milder
letter of Sobolev who was in Moscow.
The ``Savonarola'' letter was never mentioned.
Unfortunately, the official copy of the Resolution of
the Bureau of the Mathematics Division did not arrived at
Novosibirsk (the time of facsimile communication had not come
yet).
Alexandrov briefed the audience about
this Resolution. However, not without reason it is
said:``you're nobody till somebody gives you a~sheet of paper.''
V.~ A. Koptyug,
never  feeling anything positive towards Alexandrov,
moderated the discussion with reference to the unclear standpoint of
the Institute of Mathematics and the absence of the Moscow Resolution
in writing.
Of no avail were the vehement statements of the members of the Presidium
Academicians G.~K. Boreskov, S.~S. Kutateladze (1914--1986), and
A.~N. Skrinski\u\i{} who condemned the slander against Kolmogorov
and insisted on a principled reaction.
In result there was adopted a rather insipid resolution which stated that
the editorial staff of the newspaper made a serious mistake
by publication of Merzlyakov's article
``written in the style inadequate to the spirit and
aims of the newspaper.''
That was how slander had become a style in the opinion of a part of
the then leadership of the Siberian Division.

The efforts of the supporters of Kolmogorov
brought about a tactical success:
the Decree of the Presidium of the Supreme Soviet
of the USSR was signed on April 22 upon the decoration of
Academician A.~N. Kolmogorov with the Order of the October Revolution
for his great contributions to the development of
the science of mathematics and the long-term and fruitful
pedagogical activities on the occasion of
the 80 years of his birth.
It seems to me that Kolmogorov had  never become aware
of Merzlyakov's article.

Of great importance to Novosibirsk was the publication in
the issue of May 12 of the newspaper {\it Science in Siberia\/}
of an article about Kolmogorov which was written by
S.~L. Sobolev, A.~A. Borovkov, and V.~V. Yurinsky.
Their article ranked Kolmogorov as one
of the most eminent mathematicians on the twentieth century,
an outstanding teacher, an ardent patriot, and the founder
of his scientific school of a worldwide reputation
and few analogs in the history of science.
The authors particularly  emphasized the undisputable influence of Kolmogorov
on the development of
mathematics in Siberia.

This did not close the case however.
``The Special Opinion of
L.~S. Pontryagin'' was made public already on April~30.
In this article Pontryagin expressed his disagreement with the
Resolution of the Bureau of the Mathematics Division
(he was the member of the Bureau but missed the meeting on
March 25 since he was ill). He refuted the accusation against Merzlyakov
of slandering Kolmogorov and estimated the article
``generally in the positive since it summons up
citizenship which is in great demand of our scientists.''

In particular, Pontryagin wrote:

\medskip
\itemitem{}{\eightpoint\sl\obeyspaces\indent
I ascertain that the statement of Yu.~I. Merzlyakov
about Kolmogorov, even in deciphered form,
cannot be viewed as slander. It  does not ascertain
any causal relation between the failure of the reform
and the awarding of the prize. But the thought about
a causal relation can be borne in upon the reader.}

\itemitem{}{\eightpoint\sl\obeyspaces\indent
It was after this meeting already that I received  some responses to
the Yu.~I. Merzlyakov article. One of them showed disapproval
(by Academician A.~D.~Alexandrov) and three of them showed approval
(by Academician/writer Leonov; Mathematician, ScD V.~D.~Mazurov; and
the chiefs of the Mechanics and Mathematics Department of
NSU:\footnote{The abbreviation of Novosibirsk State University.}
Dean M.~M. Lavrent$'$ev and Secretary of the Party Bureau D.~E.~Zakharov).}
\medskip
\noindent
``The Special Opinion'' pinpointed a few rare facts of  public
subscription to soiling Kolmogorov's reputation.
Pontryagin's text full of the bits of an open polemic with Alexandrov raised
the question: ``Whom does A.~D. Alexandrov
defend so vehemently in his response?''.
There was little doubt that Alexandrov
would leave this question rhetorical.

Alexandrov finished his response to Pontryagin on May 28.
Confirming his view of Merzlyakov's article as
politically slanderous insinuation,  Alexandrov wrote:

\medskip
\itemitem{}{\eightpoint\sl\obeyspaces\indent
In my response to Merzlyakov's article I characterized this as
baseness and I reiterate:
$\underline{\text{this is baseness, and the
meanest
baseness at that}}$.}
\itemitem{}{\eightpoint\sl\obeyspaces\indent
Academician Pontryagin is not a young man and he knows
the intended consequences of such  baseness
in the times of the year  1937.
He could know in particular that
Nikola\u\i{} Nikolaevich Vavilov, a great Russian
scientist/biologist, died in prison since someone casted
a political slanderous innuendo about him.
Now Academician Pontryagin
supports the revival of political slanders and insinuations
and even discerns some ``citizenship'' in them.
However they were condemned by our party and people long ago.
It is the Bureau of the Mathematics Division
that revealed the genuine citizenship by repulsing
Merzlyakov's slander. The ``citizenship'' in the sense of Pontryagin
was  revealed already in his article in {\it The Communist\/}
where he spread slander against our mathematics.
Now it is revealed once again in his
``Special Opinion''  supporting baseness and slander
against not only A.~N. Kolmogorov
but also the whole school of  our scientists
which supposedly incorporates a crawling horde of the most monstrous careerists
and villains...}
\medskip
\noindent
The copies  of the March Resolutions of the Bureau of the Mathematics
Division and the Presidium of the  SDAS of the USSR
were displayed on the advertisement board of the
Institute of Mathematics of the SDAS of the USSR from  July 2~ to July~7.
So ended the crisis of  ``patriotically slanderous citizenship''
at Novosibirsk in~1983.

The above events in the history of science in Russia may be
compared only with the so-called ``Academician
N.~N. Luzin Affair.'' The pivotal distinction of the year
1986  from the year  1936 lies in the fact that
the personality of Kolmogorov
had morally united the overwhelming majority
of the Russian mathematicians who shielded
their professional community
from slander and political insinuation.

{\bf Sic transit separation.}

\head
Science at the Center of Culture
\endhead

A.~D. was a person with a universal outlook.
Through much suffering he did achieve a perfect system of views
that allowed him to analyze the general philosophical problems
and meet the challenges of contemporary life.

I had many opportunities to listen his public lectures
which always evoked a~vivid response of any audience.
I recall his brilliant talk
at the conference
``The Place of Science in the Modern Culture'' which was arranged in
Academgorodok near Novosibirsk in the end of April of 1987.

A.~D. titled his talk ``Science at the Center of Culture''
so biting a part of the audience with antipathy to
science. In my files there still reside  some records of the
main points of his talk. I insert a few of them here.

\medskip
\itemitem{}{\eightpoint\sl
We live in the age of science.}

\itemitem{}{\eightpoint\sl
\indent
False theses: science beyond culture;
science as next of kin to utopia and ideology;
science as a tool for dehumanization.}

\itemitem{}{\eightpoint\sl
\indent
This is a spite of philosophers.  A philosopher is an
unsuccessful scientist full of mania grandiose.}

\itemitem{}{\eightpoint\sl
\indent
Science occupies the center of culture. Objectively, science is
a system of knowledge and conception...
Man must stand in the center of science. Man is not only
a creator but also an object and ultimate aim
of research and thought. Science asks not only
``How?'' but also ``What for?''.}

\itemitem{}{\eightpoint\sl
\indent
Truth is a tool of good. Science leads to truth
and  its entire credenda appeal to the mind so liberating
human mentality.}
\medskip
\noindent
A.~D.~ knew much about religion, always contrasting
religious belief and scientific search.
With love to precise definitions innate in mathematicians,
he often cited the following words by
Vl.~S.~Solov$'$\"ev from the article ``Belief''
in the Encyclopedic Dictionary by
F.~A.~ Brokgauz and I.~A. Efron:
\medskip
\itemitem{}{\eightpoint\sl
\noindent
{\bf Belief} (philos.) means the acceptance of something as being true
with the resolution surpassing the power of external proofs
by fact and formal logic.}
\medskip\noindent
 A.~D. was fond of reiterating that he believes in
nothing.
This statement usually  called about the retort of the audience:
``Neither in communism?''
which always gained the affirmative answer of  A.~D.
It goes without saying that the lectures of
A.~D. were often accompanied with sneaking letters to
various local party committees.

A.~D. had explicated his views of interrelation between
religion and science in the booklet
``Scientific Search and Religious Belief''
which was published by Politizdat in 1974  with run of many thousand copies.
It seems to me that this article does not lose its actual value
nowadays in the time of an unprecedented blossom of
mysticism and pseudoscience.

\head
O.~A.~Ladyzhenskaya and a~Struggle Against the Last Insinuations
\endhead

At the end of the 1980s Aleksandr Danilovich suddenly became
a target of some slanderous attacks that ran as far as
accusations of ``lysenkoism.''
Yu.~G.~Reshetnyak and I were compelled to
write much about  A.~D.
Hatred to calumniators boiled in our soles. However,
we worked at ease feeling the inspiring
warmth of final exposition of a just-proven new theorem.
In the most critical moments of controversy
we readily found out
many objective facts witnessing the intellectual
honesty of A.~D. and his devotion to serving
science and taking care of the fates of his fellow scientists.

Stuffed up with concocted reminiscences,
massaged citations, archived data
full of sneaking and quasi-sneaking letters
to ``competent authorities''
and having mastered up many tricks typical of a barrister,
I grew up to appraise the moral standpoint
of O.~A.~Ladyzhenskaya tied with  A.~D.
by many years of friendly relations.

In the spring of 1989 I happened to be in Leningrad
at the peak of controversy about Alexandrov's
``lysenkoism.''
Olga Aleksandrovna asked me to visit her in
LOMI (the Russian abbreviation for the Leningrad Department
of the Steklov Institute).
In contrast to the majority (including some friends
of A.~D. who always requested the objective proofs
 of A.~D.'s innocence),
Olga Aleksandrovna rejected from the very beginning all my attempts to
show papers, compare figures, etc.:
``S\"ema! I need none of this stuff.
Tell me only what we must do for A.~D.''

It seemed to me that the formal position of Leningrad's mathematicians
will be important for A.~D.
Olga Aleksandrovna agreed with this opinion.
She was then a deputy chairperson of the Leningrad Mathematical Society  (LMS)
(and the chairperson was D.~K.~Faddeev).

Soon after that  V.~A. Zalgaller sent me to Novosibirsk
the following statement of the LMS which was adopted unanimously
at the meeting of March 28,~1989:

\medskip
\itemitem{}{\eightpoint\sl
Concerning the publication by the journal {\it Energy\/}
(1989:{\bf1}) of a letter of  Academician of the Siberian Division
of the Academy of Sciences of the USSR
V.~E.~Nakoryakov, the Leningrad Mathematical Society
(LMS) announces that the letter by V.~E.~Nakoryakov
contains slander (refutable by proof)
and an attempt to defame Academician A.~D. Alexandrov, a member of the LMS
and an outstanding mathematician.
The Leningrad scientists remember many good deeds by
A.~D.~Alexandrov: his efforts helped to save science and
particular scientists in the grim years which
required his great personal fortitude.}
\medskip
\noindent
A.~D. was touched with this statement. Also, it was
a great help
for
Yu.~G.~Reshetnyak and me in the public polemics of those years.

The reader seeking for more detail can restore the main events
by looking through the corresponding publications in
the {\it Herald of the Academy of Sciences of the USSR}
(1989:{\bf7}; 1990:{\bf3}) and the relevant articles
in the issues of the newspaper  {\it Science in Siberia\/}
of March 10 and October 13, 1989.

When a decade has elapsed, sharp contrast transpires
between
the figure of silence ({\it{aposiopesis}})
of the top officials of the Academy such as
V.~A.~Kirillin,  V.~A.~Koptyug, G.~I.~Marchuk, et al.
and the behavior of  the scientists who consider
the defense of the honor of a colleague against slander
as their personal duty.

I keep a few letters that were unpublished in view
of the standpoint of the then Academy  bosses.
I cherish the words of my long-term friend V. M. Tikhomirov, a professor
at Lomonosov State University in Moscow:

\medskip
\itemitem{}{\eightpoint\sl\indent
I am sure that A.~D.~Alexandrov belongs to those who have always served
the forces of good. I wish to express through your newspaper
my feeling of admiration for him,  his
brilliance, intellectual gift, and human generosity.
I've never heart that Aleksandr Danilovich
caused harm to the persons he met in life but I heard that
he helped them and promoted the development of science.}

\itemitem{}{\eightpoint\sl
\indent Of utmost importance for me are the  words by
V.~I.~Smirnov, a person of unsurpassable moral standards,
who wrote that A.~D. Alexandrov controlled the University
using the power of  moral authority!}
\medskip
\noindent
There is  no denying that the attitude of contemporaries
meant much to A.~D.

I have no desire to expatiate upon this story
even though it had a ``happy end'':
In October of 1990 A.~D., the only mathematician in a group of biologists,
was decorated for  special contribution to preservation and development
of genetics and selection in this country.

The Decoration Decree appeared by the initiative of Professor
N.~N.~Vorontsov who then hold the position of
the Chairman of the Governmental Committee for Nature of the USSR.
In a lengthy interview to the  newspaper {\it Izvestiya\/}
as of November~3, 1990 Nikola\u\i{} Nikolaevich testified:
\medskip
\itemitem{}{\eightpoint\sl\indent
Aleksandr Danilovich was the Rector of Leningrad State University
and he made much for preservation and development of genetics.
He invited to LSU many of those expelled
for their scientific views from other cities.
Young persons simply fled to Aleksandr Danilovich Alexandrov
to gain custody. The courses of lectures in LSU
differed drastically from the Lysenko rigmarole that was
delivered (and, I am afraid of that, is still delivered)
by the teachers of agricultural colleges.
This determined the general atmosphere of the academic life of Leningrad.}

\itemitem{}{\eightpoint\sl\indent
Alexandrov took care of the level of
science as a whole. All scientists know:
liquidation of one of the branches  will bring about
repercussions on the entire  frontiers of science.
That is why in many running years, many physicists and mathematicians
wrote letters to the Central Committee of the Party
about the importance of genetics. By the way, when
somebody says that A.~D.~Sakharov was late in taking the road of political
struggle,
it is not true.  His name  appeared in the letter of
physicists  of 1953  together with the names of
Kapitsa, Sem\"enov, and Varga. This letter was handled to Khrushch\"ev
by Kurchatov. The letter of physicists was followed by
a letter of mathematicians: Kolmogorov, Sobolev,
Alexandrov, and Lavrent$'$ev. I was a first-year postgraduate
when I collected their signatures.}
\bigskip

\head
The English Language
\endhead

Another not universally known trait deserves mentioning.
A.~D. was a person of a discriminating  artistic taste
with a poetic gift.
He wrote many poems and plays but most of them are lost since he
had an absolute memory and  wrote them down only
at somebody's request or to make a present of them.

A.~D. was in full command of the English language:
He delivered lectures, cited classics, and even
wrote poems in English.
S.~I. Zalgaller saved in her memory the following
lines:
\itemitem{}{\obeylines\obeyspaces\eightpoint\it
             My heart is full of burning wishes,
             My soul is under spell of thine,
             Kiss me: yours kisses are delicious
             More sweet to me than myrrth and wine.
             Oh lean against my heart with mildness,
             And I shall dream in happy silence,
             Till there will come the joyful day
             And gloom of night will fly away.}

\medskip
\noindent
Not later than in 1944 A.~D. had made this interpretation of a celebrated poem
that was written by  A.~S. Pushkin in Russian as far back as in 1825
and soon became an immortal romance  by M.~I. Glinka.

It is curious but one of our first conversations in the mid-1960s
ran in English
(I was a freshmen; and  ~A.~D., a brand-new academician).
I recall the presence of some ``English-speaking'' diplomat
in the hall of a small canteen in the Golden Valley
where we dined
those years.  A.~D. remarked that it is indecent to
use the language that is not comprehended by everyone
present and we proceeded in English.

I also recall an episode of the 1970s when on some occasion  I cited
a few lines of the 66th sonnet by W.~Shakespeare in English,
and  A.~D. continued recital in a flash.
This took place long before the famous Georgian  ``Repentance'' by T.~ Abuladze.

The circumstances of the beginning of the 1990s
drove me to writing a short booklet on
English grammar to alleviate the burdens of life
for my friends who were seeking  some sources of nourishment.
A.~D., always a very attentive reader,
pinpointed a slip in a~King James's citation of
Ecclesiastes.

And in June of 1993 he sent me
the following lines in a sloppy handwriting:
\itemitem{}{\obeylines\obeyspaces\eightpoint\it
Since legs, nor hands, nor eyes, nor strong creative brain,
But weakness and decay oversway their power,
I am compelled forever to refrain
From everything but waiting for my hour.
}
\medskip
\noindent
He has never sent me any verses since then...


\enddocument